\documentclass[12pt, a4paper]{article}
\usepackage[top=1.2in, bottom=1.2in, left=0.9in, right=0.9in]{geometry}
\usepackage[small]{titlesec}
\titleformat{\subsubsection}[runin]{\bfseries}{\upshape\thesubsubsection}{0.5em}{} 

\usepackage{amssymb, amsmath, amsfonts}
\usepackage{mathabx, pdfpages}
\usepackage{tikz}
\usepackage{tikz-cd}

\usepackage{hyperref} 
\hypersetup{
	colorlinks=true,
	linkcolor=red,
	citecolor=blue,
	filecolor=cyan,      
	urlcolor=magenta}

\usepackage{amsthm}
\newtheoremstyle{remark}{}{}{\upshape}{}{\bfseries}{.}{0.5em}{}
\newtheorem{theorem}{Theorem}[section]

\newtheorem{corollary}[theorem]{Corollary}
\newtheorem{proposition}[theorem]{Proposition}
\newtheorem{definition}[theorem]{Definition}
\theoremstyle{remark}

\newtheorem{remark}[theorem]{Remark}

\newtheorem{para}[theorem]{ } 

\usepackage[all,cmtip]{xy}
\usepackage{bm}
\usepackage[utf8]{inputenc}

\def\id{{\rm id}}
\long\def\nodo#1{{}}
\def\op{\mathrm{op}}
\def\gg{\mathfrak{g}}

\def\MR#1{} 

\def\btr{\blacktriangleright}
\def\btl{\blacktriangleleft}

\DeclareMathOperator{\hotimes}{ \hat\otimes }
\DeclareMathOperator{\totimes}{ \tilde\otimes }

\def\V{\mathrm{Vect}}

\def\pV{\mathrm{pro}\V}

\def\iV{\mathrm{ind}\V}

\def\ipV{\mathrm{ind}\mathrm{pro}\V}

\def\Ug{U(\mathfrak{g})}

\def\VV{\mathcal{V}}

\begin{document}
\title{\bf Internal Hopf algebroid}

\author{ Martina Stoji\'c  }

\date{ \normalsize \today}
\maketitle

\begin{abstract} We introduce a natural generalization of the definition of a symmetric Hopf algebroid, internal to any symmetric monoidal category with coequalizers that commute with the monoidal product. Motivation for this is the study of Heisenberg doubles of countably dimensional Hopf algebras $A$ as internal Hopf algebroids over a (noncommutative) base $A$ in the category $\mathrm{indproVect}$ of filtered cofiltered vector spaces introduced by the author. One example of such Heisenberg double is internal Hopf algebroid  $U(\mathfrak{g}) \sharp U(\mathfrak{g})^*$  over universal enveloping algebra $U(\mathfrak{g})$ of a finite-dimesional Lie algebra $\mathfrak{g}$ that is a properly internalized version of a completed Hopf algebroid previously studied as a Lie algebra type noncommutative phase space.
\end{abstract}

\tableofcontents

\section{Introduction}

This article is about internalization of notions of bialgebroid and Hopf algebroid in any symmetric monoidal category with coequalizers that commute with the monoidal product. This alowed the author to study Heisenberg doubles of countably dimensional Hopf algebras~$A$  as internal Hopf algebroids over a (noncommutative) base $A$ in the category $\ipV$ of filtered cofiltered vector spaces introduced in author's dissertation \cite{CHA}. One example of such Heisenberg double is internal Hopf algebroid  $U(\mathfrak{g}) \sharp U(\mathfrak{g})^*$  over universal enveloping algebra $U(\mathfrak{g})$ of a finite-dimesional Lie algebra $\mathfrak{g}$ that is a properly internalized version of a completed Hopf algebroid previously studied as a Lie algebra type noncommutative phase space in \cite{halgoid}. The axioms in \cite{halgoid} deal with the problem of completions by introducing \emph{ad hoc} versions of axioms that are not properly internalized to some monoidal category. 
 
Commutative Hopf algebroids arise as objects dual to groupoids. Namely, cogroupoid in the symmetric monoidal category of commutative $k$-algebras is commutative Hopf algebroid over a commutative base algebra. Commutative bialgebroid lacks only one structure morphism, antipode, to be a commutative Hopf algebroid.

Noncommutative generalization of commutative bialgebroids are associative $k$-bi\-al\-ge\-bro\-ids from \cite{Lu}, whose axioms were further studied and discussed in \cite{xu,BrzMilitaru,bohmHbk,SS:two,Stojic}. In Section \ref{sec:internalbio}, we exhibit in detail the internalized version of that notion due to G.~B\"ohm~\cite{bohmInternal} in the context of a symmetric monoidal category with coequalizers that commute with the monoidal product. On this internalization of bialgebroid, we base our definition of internal Hopf algebroid in Section \ref{sec:internalHA}.  Justification of the definition itself requires more care than in the category of vector spaces, where the definition of symmetric Hopf algebroid is developed by G.~B\"ohm~\cite{Bohmalt,Bohm}, correcting some shortcomings of the unsymmetrical Hopf algebroids according to J-H.~Lu~\cite{Lu}.

In Appendix \ref{summary}, the summary of author's dissertation is given where the main applications of the results are. This material will appear elsewhere. 

\section{Internal bialgebroid}\label{sec:internalbio}
In this article $\mathcal{V} = (\mathcal{V}, \otimes, k, a, l, r, \tau)$ is a fixed symmetric monoidal category with coequalizers that commute with the monoidal product $\otimes$ in the sense of Definition~\ref{def:koujedkom}. These assumptions are essential for the definition of internal bialgebroid.

\begin{definition}\label{def:koujedkom}
	Let $\mathcal{V}=(\mathcal{V},\otimes,k,a,l,r)$ be a monoidal category with coequalizers, that is, in which every pair of morphisms has a coequalizer. We say that \emph{coequalizers in $\mathcal{V}$ commute with the tensor product} $\otimes$ if for every for every parallel pair $f,g \colon X\to Y$ of morphisms in $\VV$ the following holds: 
	if $\xi \colon Y\to Z$ is a coequalizer of the pair $f,g$
	$$\begin{tikzcd}
	X
	\ar[shift left=2]{r}{f} \ar[swap,shift right=1]{r}{g} &
	Y \ar{r}{\xi}& Z,
	\end{tikzcd}$$
	then for every object $W$ in category $\mathcal{V}$ the morphism $\xi\otimes\id_W$ is a coequalizer of the pair $f\otimes\id_W, g\otimes\id_W\colon X\otimes W\to Y\otimes W$ 
	$$\begin{tikzcd}
	X\otimes W
	\ar[shift left=2]{r}{f\otimes\id_W} \ar[shift right=1, swap]{r}{g\otimes\id_W}&
	Y\otimes W \ar{r}{\xi \otimes \id_W}& Z\otimes W,
	\end{tikzcd}$$
	and the morphism $\id_W\otimes\xi$ is a coequalizer of the pair $\id_W\otimes f,\id_W\otimes g\colon W\otimes X\to W\otimes Y$, 
	$$\begin{tikzcd}
	W\otimes X
	\ar[shift left=2]{r}{\id_W \otimes f} \ar[shift right=1, swap]{r}{\id_W \otimes g}&
	W\otimes Y \ar{r}{\id_W \otimes \xi}& W\otimes Z.
	\end{tikzcd}$$
\end{definition}
We use the fact that the tensor product of monoids in symmetric monoidal category is again a monoid in that category.
\begin{proposition} \label{prop:tenzmonoida}
	Let $(R,\mu_R,\eta_R)$ and $(S,\mu_S,\eta_S)$ be monoids in symmetric monoidal category $\mathcal{V} = (\mathcal{V},\otimes,k,a,l,r,\tau)$. Then $(R\otimes S,\mu_{R\otimes S},\eta_{R\otimes S})$, where
	$$
	\mu_{R\otimes S} = (\mu_R\otimes\mu_S)\circ(\id_R\otimes\tau_{S,R}\otimes\id_S),
	\quad \eta_{R\otimes S} = (\eta_R\otimes\eta_S)\circ l_k^{-1},
	$$
	is again a monoid in $\mathcal{V}$.
\end{proposition}
\begin{remark} It is the same if we choose to use the left or the right unit in the previous proposition, because in every monoidal category we have $l_k = r_k \colon k\otimes k \to k$. The proof is standard, see ~\cite{Kelly}. In abstract Sweedler notation the multiplication and the unit from the previous proposition are
$$(r\otimes s )\cdot (r'\otimes s') = rr'\otimes ss', \quad 1_{R\otimes S} = 1_R \otimes 1_S,$$
and therefore this multiplication is often called componentwise multiplication.
\end{remark}
We present now the definition of G.~B\"ohm of internal bialgebroid in such categories \cite{bohmInternal} in detail. 
\subsection{Preparatory results}\label{ss:unutbialglem}
In this subsection, each proposition uses the notation introduced in previous propositions and between propositions, and relies on properties of denoted objects that are proven in the previous propositions.

\begin{para} Let $(L,\mu_L,\eta_L)$ i $(H,\mu_H,\eta_H)$ be monoids in $\mathcal{V}$. Since in $\VV$ coequalizers commute with the monoidal product, the monoidal category of $L$-bimodules $({}_L \mathcal{V}{}_L, \otimes_L, L)$ is well defined.
\end{para}

\begin{proposition}\label{prop:alphabeta}
	Let $\alpha\colon L \to H$ and $\beta \colon L^{\op} \to H$ be morphisms of monoids in $\mathcal{V}$ such that 
	\begin{equation}\label{eq:commalphabeta}
	\mu_H \circ \tau_{H,H} \circ (\alpha \otimes \beta) = \mu_H \circ (\alpha \otimes \beta).
	\end{equation}
	Then $H$ is an internal $L$-bimodule by left action $\mu_H \circ (\alpha \otimes \id_H)$ and right action $\mu_H \circ \tau_{H,H} \circ (\id_H \otimes \beta)$. 
\end{proposition}
\begin{proof}
	We first prove that $\mu_H \circ (\alpha \otimes \id_H)$ is a left action. We have that
	$$\begin{array}{lcl}
	\mu_H\circ(\alpha\otimes(\mu_H\circ(\alpha\otimes\id_H)))
	&=& 
	\mu_H\circ((\mu_H\circ(\alpha\otimes\alpha))\otimes\id_H)\\
	&=& \mu_H\circ((\alpha\circ\mu_L)\otimes\id_H),
	\end{array}$$
	where in the first line we used associativity of $\mu_H$, and in the second line that $\alpha$ is a morphism of monoids. Further, for proof of the compatibility with unit $\eta_H$, simply notice that $\mu_H\circ((\alpha\circ\eta_L)\otimes\id_H) = 
	\mu_H\circ(\eta_H\otimes\id_H) = l_H$. 
	Claim that $\mu_H \circ \tau_{H,H} \circ (\id_H \otimes \beta)$ is a right action is proven analogously.
	
	Finally, we have to check that the left and the right action commute, that is,
	$$
	\mu_H\circ(\alpha\otimes(\mu_H\circ\tau_{H,H}\circ(\id_H\otimes\beta)))
	= \mu_H\circ\tau_{H,H}\circ((\mu_H\circ(\alpha\otimes\id_H))\otimes
	\beta) 
	$$
	as morphisms $L\otimes H\otimes L\to H$. Left hand side is equal to $\mu_H\circ(\mu_H\otimes\id_H)\circ(\alpha\otimes\beta\otimes\id_H)\circ(\id_L\otimes\tau_{H,L})$. Right hand side is equal to
	$\mu_H\circ(\beta\otimes(\mu_H\circ(\alpha\otimes\id_H)))\circ\tau_{L\otimes H,L},$ which is $\mu_H\circ((\mu_H \circ (\beta\otimes\alpha))\otimes\id_H)\circ\tau_{L\otimes H,L}$. By the assumption~(\ref{eq:commalphabeta}), this is equal to	$\mu_H\circ((\mu_H\circ(\alpha\otimes \beta)\circ\tau_{L,L})\otimes\id_H)\circ
	\tau_{L\otimes H,L}$, which is $ \mu_H\circ(\mu_H\otimes \id_H) \circ (\alpha\otimes \beta \otimes \id_H ) \circ (\tau_{L,L}\otimes\id_H)\circ
	\tau_{L\otimes H,L}$. The equality now follows from
	$(\tau_{L,L}\otimes\id_H)\circ\tau_{L\otimes H,L} = \id_L\otimes\tau_{H,L}$.
\end{proof}

\begin{para}
	For $\alpha$ and $\beta$ as in Proposition \ref{prop:alphabeta}, denote by $\tilde\alpha$ and $\tilde\beta$ the compositions:
	$$
	\begin{tikzcd}
	L \ar{r}{\cong} \ar[bend right]{rr}{\tilde\alpha} & k\otimes L \ar{r}{\eta_H \otimes \alpha}  & H \otimes H  & L \ar{r}{\cong} \ar[bend right]{rr}{\tilde\beta}  & L \otimes k \ar{r}{\beta \otimes \eta_H } & H \otimes H.
	\end{tikzcd}
	$$ 
	In abstract Sweedler notation, this is: $\tilde{\alpha}(l)= 1_H  \otimes \alpha(l)$ i $\tilde{\beta}(l) = \beta(l) \otimes 1_H$.  
\end{para}
\begin{para}
	Denote $\mu_{H\otimes H} := (\mu_H \otimes \mu_H)\circ (\id_H \otimes \tau_{H,H}  \otimes \id_H)$ and  $\eta_{H\otimes H} := (\eta_H \otimes \eta_H) \circ l_k^{-1}$. Then $(H\otimes H,\mu_{H\otimes H},\eta_{H\otimes H})$ is a monoid in category $\mathcal{V}$ by Proposition~\ref{prop:tenzmonoida}.
\end{para}
\begin{para}
	Let $\pi$ be the canonical morphism $H\otimes H\to H\otimes_L H$, that is, the coequalizer 
	$$
	\begin{tikzcd}
		H \otimes L \otimes H \ar[shift left=2]{r} \ar[shift right]{r} & H \otimes H  \ar{r}{\pi} & H \otimes_L H 
	\end{tikzcd}
	$$
	for $L$-bimodule $H$. It is easy to see that then $\pi$ is also the coequalizer of the parallel pair of morphisms $\mu_{H\otimes H}\circ(\tilde\alpha\otimes\id_{H\otimes H})$, $\mu_{H\otimes H}\circ(\tilde\beta\otimes\id_{H\otimes H})$,
	$$
	\begin{tikzcd}
		L \otimes (H \otimes H) \ar[shift left=2]{r}{\tilde\alpha\otimes\id} \ar[shift right, swap]{r}{\tilde\beta\otimes\id} & (H\otimes H)\otimes (H\otimes H) \ar{r}{\mu_{H\otimes H}} & H \otimes H  \ar{r}{\pi} & H \otimes_L H.
	\end{tikzcd}
	$$
\end{para}

\begin{proposition} 
	(i) There exists a unique $\rho $ such that the following diagram is commutative:
	$$
	\begin{tikzcd}
		(H \otimes H) \otimes (H \otimes H)  \ar{r}{\mu_{H\otimes H}} 
		\ar{d}{ \pi \otimes \id}& 
		(H\otimes H) 
		\ar{d}{ \pi}\\
		(H \otimes_L H) \otimes (H \otimes H)  \ar[dashed]{r}{\rho} 
		& 
		(H\otimes_L H) .
	\end{tikzcd}
	$$
	(ii) This unique $\rho$ is a right action.
\end{proposition}
\begin{proof}
	The uniqueness of $\rho$ follows from the following diagram.
	$$
	\begin{tikzcd}
		L \otimes (H \otimes H) \otimes (H \otimes H) 
		\ar{r}{\id\otimes \mu_{H\otimes H}} \ar[shift right=2, swap]{d}{\tilde\alpha\otimes\id\otimes\id} \ar[shift left]{d}{\tilde\beta\otimes\id\otimes \id} & 
		L \otimes (H\otimes H) 
		\ar[shift right=2, swap]{d}{\tilde\alpha\otimes\id} 
		\ar[shift left]{d}{\tilde\beta\otimes\id}  
		\\
		(H \otimes H) \otimes (H \otimes H) \otimes (H \otimes H) \ar{r}{\id\otimes \mu_{H\otimes H}} 
		\ar{d}{ \mu_{H\otimes H} \otimes \id}& 
		(H \otimes H) \otimes (H\otimes H) 
		\ar{d}{ \mu_{H\otimes H}}
		\\
		(H \otimes H) \otimes (H \otimes H)  \ar{r}{\mu_{H\otimes H}} 
		\ar{d}{ \pi \otimes \id}& 
		(H\otimes H) 
		\ar{d}{ \pi}\\
		(H \otimes_L H) \otimes (H \otimes H)  \ar[dashed]{r}{\rho} 
		&
		(H\otimes_L H) 
	\end{tikzcd}
	$$
	Vertical maps on the right comprise a diagram for coequalizer $\pi$, and vertical maps on the left a diagram for coequalizer $\pi\otimes \id$, because in the category $\VV$ coequalizers commute with the monoidal product. Top rectangle diagram is sequentially commutative because of functoriality of the monoidal product, and middle rectangle diagram because $\mu_{H\otimes H}$ is associative. Therefore, we conclude that $\pi \circ \mu_{H\otimes H}$ coequalizes the parallel pair of vertical morphisms on the left. From this  it follows that there exists a unique $\rho$, because of the universal property of coequalizer $\pi \otimes \id$. 
	
	We now prove that $\rho$ is a right action.
	$$
	\begin{tikzcd}
		(H\!\otimes_L\! H) \otimes (H\!\otimes\! H) \otimes (H\!\otimes\! H) 
		\ar{r}{\rho\otimes\id} \ar[bend right=80, swap]{ddd}{\id\otimes \mu_{H\otimes H}} & 
		(H\!\otimes_L\! H) \otimes (H\!\otimes\! H) 
		\ar[bend left=80]{ddd}{\rho} \\
		(H\!\otimes\! H) \otimes (H\!\otimes\! H) \otimes (H\!\otimes\! H) 
		\ar{r}{\mu_{H\otimes H}\otimes\id} \ar{d}{\id\otimes \mu_{H\otimes H}} \ar{u}{\pi\otimes\id\otimes\id} &
		(H\!\otimes\! H) \otimes (H\!\otimes\! H) 
		\ar{d}{\mu_{H\otimes H}} \ar{u}{\pi\otimes\id}\\
		(H\!\otimes\! H) \otimes (H\!\otimes\! H) 
		\ar{r}{\mu_{H\otimes H}} \ar{d}{\pi\otimes\id} &
		(H\!\otimes\! H)  
		\ar{d}{\pi}\\
		(H\!\otimes_L\! H) \otimes (H\!\otimes\! H) 
		\ar{r}  & 
		(H\!\otimes_L\! H)
	\end{tikzcd}
	$$
	Commutativity of five internal quadrangle diagrams follows from the associativity of $\mu_{H\otimes H}$, definition of $\rho$ and functoriality of the monoidal product. Since $\pi\otimes\id\otimes\id$ is an epimorphism, we conclude that the outer quadrangle diagram is also commutative.
	
	$$
	\begin{tikzcd}
		(H\!\otimes\! H) \otimes k \ar[bend left]{rr} \ar{r}{\id\otimes\eta_{H\otimes H}} \ar{d}{\pi\otimes\id} & (H\!\otimes\! H) \otimes (H\!\otimes\! H) \ar{r}{\mu_{H\otimes H}} \ar{d}{\pi\otimes\id} & (H\!\otimes\! H) \ar{d}{\pi}\\
		(H\!\otimes_L\! H) \otimes k \ar{r}{\id\otimes\eta_{H\otimes H}} \ar[bend right]{rr} & (H\!\otimes_L\! H) \otimes (H\!\otimes\! H) \ar{r}{\rho} & (H\!\otimes_L\! H) \\
		&&
	\end{tikzcd}
	$$
	Commutativity of the top triangle, inner quadrangles and outer quadrangle diagram follow from the compatibility of unit and counit in monoid $H \otimes H$, functoriality of the monoidal product, definition of $\rho$ and the fact that the right unitor in a category is a natural transformation. Since $\pi\otimes\id$ is an epimorphism, we conclude that the bottom triangle diagram is commutative.
\end{proof}

\begin{para} Further, let $_L H_L$ be equipped with the structure of comonoid $(H, \Delta, \epsilon)$ in the monoidal category of $L$-bimodules. That is, $\Delta \colon H \to H \otimes_L H$ is a coassociative $L$-bimodule morphism and  $\epsilon \colon H \to L$ is a counit for $\Delta$, also an $L$-bimodule morphism.
\end{para}
\begin{proposition}\label{prop:rholambda} 
	Assume that $\rho$ coequalizes the pair of morphisms $\Delta\otimes \tilde\alpha$ and $\Delta\otimes\tilde\beta$, 
	$$
	\begin{tikzcd}
		H \otimes L  \ar{d}{\Delta\otimes \id}  & & \\
		(H\! \otimes_L\! H) \otimes L 
		\ar[shift left=2]{r}{\id\otimes\tilde\alpha} \ar[shift right, swap]{r}{\id\otimes\tilde\beta} 
		& 
		(H\! \otimes_L\! H) \otimes (H\! \otimes\! H)  \ar{d}{\rho} &   \\
		& (H\!\otimes_L\! H). &
	\end{tikzcd}
	$$
	Then there exists a unique $\lambda$ such that the following  diagram is commutative.
	$$
	\begin{tikzcd}
		H \otimes (H\!\otimes\! H) \ar{r}{\id\otimes\pi} \ar{d}{\Delta\otimes\id} & H \otimes (H\!\otimes_L\! H) \ar[dashed]{d}{\lambda} \\
		(H\!\otimes_L\! H) \otimes (H\!\otimes\! H) \ar{r}{\rho} & (H\!\otimes_L\! H)
	\end{tikzcd}
	$$
\end{proposition}
\begin{proof}
	$$
	\begin{tikzcd}
		& &  H\otimes (H\!\otimes_L\!  H) \ar[dashed, bend left=70]{ddd}{\lambda} \\
		H \otimes L \otimes (H\! \otimes\! H) 
		\ar{d}{\Delta\otimes \id\otimes\id} \ar[shift left=2]{r}{\id\otimes\tilde\alpha\otimes\id} \ar[shift right, swap]{r}{\id\otimes\tilde\beta\otimes\id} & 
		H \otimes (H\! \otimes\! H) \otimes (H\! \otimes\! H) \ar{r}{\id\otimes \mu_{H\otimes H}} \ar{d}{\Delta\otimes \id\otimes\id} & H\otimes (H\!\otimes\! H) \ar{u}{\id\otimes\pi} \ar{d}{\Delta\otimes \id}
		\\
		(H\! \otimes_L\! H) \otimes L \otimes (H\! \otimes\! H) 
		\ar[shift left=2]{r}{\id\otimes\tilde\alpha\otimes\id} \ar[shift right, swap]{r}{\id\otimes\tilde\beta\otimes\id} 
		& 
		(H\! \otimes_L\! H) \otimes (H\! \otimes\! H) \otimes (H\! \otimes\! H) \ar{r}{\id\otimes \mu_{H\otimes H}} \ar{d}{\rho\otimes \id}& (H\! \otimes_L\! H) \otimes (H\! \otimes\! H) \ar{d}{\rho}  \\
		& (H\!\otimes_L\! H) \otimes (H\!\otimes\! H) \ar{r}{\rho} & (H\!\otimes_L\! H) 
	\end{tikzcd}
	$$
	Since in category $\VV$ coequalizers commute with the monoidal product, $\id \otimes \pi$ is a coequalizer of the morphisms on the top border of the diagram. From the assumption of the proposition, it follows that $\rho \otimes \id$ coequalizes pair $\Delta \otimes \tilde{\alpha} \otimes \id$ and $\Delta\otimes \tilde{\beta} \otimes \id$, which is on the diagram represented by the stepwise left border of the diagram. Three rectangle diagrams are commutative because of functoriality of $\otimes$ and because $\rho$ is an action (top left sequentially). It follows that $\rho \circ (\Delta \otimes \id)$ coequalizes horizontal morphisms on the top border of the diagram. Because of the universal property of coequalizer $\pi \otimes \id$, we conclude that there exists a unique morphism~$\lambda$.
\end{proof}
\begin{proposition}
	Morphisms $\lambda \otimes \id$ and $\id\otimes\rho$ commute.
\end{proposition}
\begin{proof}
	$$
	\begin{tikzcd}
		H \otimes (H\!\otimes\! H) \otimes (H\!\otimes\! H) 
		\ar{r}{\Delta\otimes\id\otimes\id} \ar{d}{\id\otimes\pi\otimes\id}
		\ar[bend right=80]{ddd}{\id\otimes \mu_{H\otimes H}} &
		(H\!\otimes_L\! H) \otimes (H\!\otimes\! H) \otimes (H\!\otimes\! H)
		\ar{d}{\rho\otimes \id}
		\ar[bend left=75]{ddd}{\rho\otimes \mu_{H\otimes H}} \\
		H \otimes (H\!\otimes_L\! H) \otimes (H\!\otimes\! H) 
		\ar{r}{\lambda\otimes\id} \ar{d}{\id\otimes\rho} &
		(H\!\otimes_L\! H) \otimes (H\!\otimes\! H) 
		\ar{d}{\rho}\\
		H \otimes (H\!\otimes_L\! H)  
		\ar{r}{\lambda}  &
		(H\!\otimes_L\! H) \\
		H \otimes (H\!\otimes\! H)  
		\ar{r}{\Delta\otimes\id} \ar{u}{\id\otimes\pi} &
		(H\!\otimes_L\! H) \otimes (H\!\otimes\! H) 
		\ar{u}{\rho}
	\end{tikzcd}
	$$
	The following subdiagrams are commutative: top rectangle (definition of $\lambda$, tensoring on the right with $\id$), quadrangle on the right ($\rho$ is an action), outer quadrangle (functoriality of the monoidal product), bottom rectangle (definition of $\lambda$) and quadrangle on the left (definition of $\rho$, tensoring on the left with $\id$) and $\id\otimes\pi\otimes\id$ is an epimorphism. Therefore, the central rectangle diagram is also commutative, proving the commutativity of $\lambda \otimes \id$ and $\id\otimes \rho$.
\end{proof}

\begin{proposition}
	If  
	\begin{equation}\label{eq:Delta11}
	\Delta\circ\eta_H = \pi\circ \eta_{H\otimes H},
	\end{equation}
	that is, the following diagram is commutative 
	$$
	\begin{tikzcd}
		k \ar{r}{\eta_H} \ar{d}{\eta_{H\otimes H}} &  H \ar{d}{\Delta} \\
		(H\!\otimes\! H)  \ar{r}{\pi} & (H\!\otimes_L\! H),
	\end{tikzcd}
	$$
	and morphism $\lambda$ satisfies
	\begin{equation}\label{eq:lambdauvjet}
	\Delta\circ \mu_H = \lambda\circ (\id_H\otimes\Delta),
	\end{equation}
	that is, the following diagram is commutative
	$$
	\begin{tikzcd}	
		H\otimes H \ar{r}{\mu_H} \ar{d}{\id\otimes\Delta} &  H \ar{d}{\Delta} \\
		H\otimes (H\!\otimes_L\! H) \ar{r}{\lambda} & (H\!\otimes_L\! H),
	\end{tikzcd}
	$$
	then $\lambda$ is a left $H$-action on $H\otimes_L H$.
\end{proposition}
In abstract Sweedler notation these two conditions are written as: $\Delta(1_H) = 1_H\otimes_L 1_H$ and $\lambda(h \otimes h'_{(1)} \otimes h'_{(2)}) = \Delta(hh')$, for $h,h'\in H$. In the category of vector spaces the second condition is simply 
$\Delta(h h') = h_{(1)} h'_{(1)}\otimes_L h_{(2)}h'_{(2)}$ for all $h,h'\in H$, but for the right side to be well defined the above axioms are required.
\begin{proof}
	\begin{equation}\label{eq:lambdaproof1}
	\begin{tikzcd}
		H \otimes H \otimes (H\!\otimes\! H)
		\ar{r}{\id\otimes\Delta\otimes\id} \ar{d}{\id\otimes\id\otimes\pi} 
		\ar[bend right=75]{ddd}{ \mu_{H\otimes H}\otimes \id} &
		H \otimes (H\!\otimes_L\! H) \otimes (H\!\otimes\! H)
		\ar{d}{\id\otimes\rho} \ar[bend left=75]{ddd}{\lambda \otimes \id} \\
		H \otimes H \otimes (H\!\otimes_L\! H)
		\ar{r}{\id\otimes\lambda} \ar{d}{\mu_{H\otimes H}\otimes\id} &
		H \otimes (H\!\otimes_L\! H)
		\ar{d}{\lambda}\\
		H \otimes (H\!\otimes_L\! H)
		\ar{r}{\lambda}  &
		(H\!\otimes_L\! H) \\
		H \otimes H \otimes (H\!\otimes\! H)
		\ar{r}{\Delta\otimes\id} \ar{u}{\id\otimes\pi} &
		(H\!\otimes_L\! H) \otimes (H\!\otimes\! H)
		\ar{u}{\rho}
	\end{tikzcd}
	\end{equation}
	The following subdiagrams are commutative: top rectangle (definition of $\lambda$, tensoring on the left with $\id$), bottom rectangle (definition of $\lambda$), quadrangle on the left (functoriality of the monoidal product), quadrangle on the right (commutativity of $\lambda\otimes \id$ and $\id\otimes \rho$) and outer quadrangle (the condition satisfied by $\lambda$ and $\Delta$ tensored on the right with $\id$), and $\id\otimes\id\otimes\pi$ is an epimorphism, hence the central rectangle diagram is also commutative, which proves the action associativity axiom for~$\lambda$.
	\vspace{-2em}
	\begin{equation}\label{eq:lambdaproof2}
	\begin{tikzcd}
		& (H\!\otimes\! H) \otimes (H\!\otimes\! H) 
		\ar{dr}{\pi\otimes\id} \ar{r}{\mu_{H\otimes H}} & 
		(H\!\otimes\! H) 
		\ar[bend left=100]{dd}{\pi} \\
		k\otimes (H\!\otimes\! H) \ar[bend left=60]{urr} \ar{ru}{\eta_{H\otimes H}\otimes \id} \ar{r}{\eta_H\otimes\id} \ar{d}{\id\otimes\pi} & H\otimes (H\!\otimes\! H) \ar{r}{\Delta\otimes\id} \ar{d}{\id\otimes\pi} & (H\!\otimes_L\! H) \otimes (H\!\otimes\! H) \ar{d}{\rho} \\
		k\otimes (H\!\otimes_L\! H) \ar[bend right=15]{rr} \ar{r}{\eta_H\otimes\id}  & H\otimes (H\!\otimes_L\! H) \ar{r}{\lambda}  & (H\!\otimes_L\! H)
	\end{tikzcd}
	\end{equation}
	The following subdiagrams are commutative: outer quadrangle (naturality of the left unit of the monoidal category), top triangle (axiom of the unit for multiplication $\mu_{H\otimes H}$), inner upper quadrangle (the condition that is satisfied by $\Delta $ and $\eta_{H}$ tensored on the right with $\id$), inner lower left rectangle (functoriality of $\otimes$), inner lower right rectangle (definition of $\lambda$), inner right quadrangle (associativity axiom for action $\rho$), and $\id\otimes\pi$ is an epimorphism, hence the bottom triangle is also commutative, which proves the action unit axiom for $\lambda$.
\end{proof}
\begin{proposition}
	Conversely, if $\lambda$ is a left action, then   $$\Delta\circ\eta_H = \pi\circ \eta_{H\otimes H},$$ $$\Delta\circ \mu_H = \lambda\circ (\id_H\otimes\Delta).$$
\end{proposition}
\begin{proof}
	We consider again the diagram~(\ref{eq:lambdaproof1}). If the central rectangle diagram is commutative, from the commutativity of other four inner quadrangle diagrams it follows that the outer quadrangle diagram is commutative, which is a diagram for the property for $\Delta$ and $\lambda$ tensored with $\id$. From this, by using the unit $\eta_{H\otimes H}$ and the fact that it is a monomorphism, we can get the commutative diagram that represents the property $\Delta\circ \mu_H = \lambda\circ (\id_H\otimes\Delta)$. 
	
	We consider again diagram~(\ref{eq:lambdaproof2}). If the bottom triangle diagram is commutative, then the inner upper quadrangle diagram is commutative, and this is property for $\Delta$ and $\pi$ tensored with $\id$. From this, by using the unit $\eta_{H\otimes H}$ and the fact that it is a monomorphism, we can get the commutative diagram representing property $\Delta\circ\eta_H = \pi\circ \eta_{H\otimes H}$. 
\end{proof}
Therefore, it is shown that the property $\Delta\circ\eta_H = \pi\circ \eta_{H\otimes H}$ is equivalent to the action unit axiom for $\lambda$, and that the property $\Delta\circ \mu_H = \lambda\circ (\id_H\otimes\Delta)$ is equivalent to the action associativity axiom for $\lambda$.
\begin{corollary} 
	$H\otimes_L H$ is an internal $H$-$(H\otimes H)$-bimodule with regard to actions $\lambda$ and~$\rho$.
\end{corollary}

\subsection{Definition of internal bialgebroid}

In the following two definitions, category $\VV$ is a symmetric monoidal category $(\mathcal{V}, \otimes, k, a, l, r, \tau)$ with coequalizers that commute with the monoidal product.

\begin{definition} \label{def:internalB} Let $L = (L, \mu_L, \eta_L)$ be a monoid in $\mathcal{V}$. \emph{Internal left $L$-bialgebroid} in category $\mathcal{V}$ is given by the following data: 
	\begin{itemize}
		\item[(i)] monoid $(H , \mu_H, \eta_H)$ in $\mathcal{V}$
		\item[(ii)] monoid morphisms  $\alpha\colon L \to H$, $\beta \colon L^{\op} \to H$ which satisfy $$\mu_H \circ \tau_{H,H} \circ (\alpha \otimes \beta) = \mu_H \circ (\alpha \otimes \beta)$$
		\item[(iii)] comonoid $(H, \Delta, \epsilon)$ in the monoidal category of $L$-bimodules $({}_L \mathcal{V}{}_L, \otimes_L, L)$, where $H$ is an internal  $L$-bimodule with regard to left action $\mu_H \circ (\alpha \otimes \id_H)$ and right action $\mu_H \circ \tau_{H,H} \circ (\id_H \otimes \beta)$
	\end{itemize}
	satisfying the following conditions:
	\begin{itemize}
		\item[(i)] Unique right action $\rho$ for which the following diagram is commutative:
		$$
		\begin{tikzcd}
			(H \otimes H) \otimes (H \otimes H)  \ar{r}{\mu_{H\otimes H}} 
			\ar{d}{ \pi \otimes \id} & 
			(H\otimes H) 
			\ar{d}{ \pi} \\
			(H \otimes_L H) \otimes (H \otimes H)  \ar[dashed]{r}{\rho} 
			& 
			(H\otimes_L H) 
		\end{tikzcd}
		$$
		satisfies
		$$\rho\circ (\Delta\otimes\beta\otimes\eta_H)\circ(\id_H\otimes r_L^{-1}) = \rho\circ (\Delta\otimes\eta_H\otimes\alpha)\circ(\id_H\otimes l_L^{-1}),$$ 
		where $r_L$ and $l_L$ are components of right unit and left unit of the monoidal category~$\mathcal{V}$, $\mu_{H\otimes H}$ is the multiplication in monoid $(H \otimes H, \mu_{H \otimes H}, \eta_{H\otimes H})$ induced by $\mu_H$ componentwise, and $\pi$ is the morphism in coequalizer $H\otimes H \to H\!\otimes_L\! H$.
		\item[(ii)] Unique morphism $\lambda$ (compare with Proposition~\ref{prop:rholambda})
		for which the following diagram is commutative: 
		$$
		\begin{tikzcd}
			H \otimes (H\!\otimes\! H) \ar{r}{\id\otimes\pi} \ar{d}{\Delta\otimes\id} & H \otimes (H\!\otimes_L\! H) \ar[dashed]{d}{\lambda} \\
			(H\!\otimes_L\! H) \otimes (H\!\otimes\! H) \ar{r}{\rho} & (H\!\otimes_L\! H)
		\end{tikzcd}
		$$
		is a left action. Equivalently, $$\Delta\circ\eta_H = \pi\circ \eta_{H\otimes H},$$ $$\Delta\circ \mu_H = \lambda\circ (\id_H\otimes\Delta)$$ hold, where $\eta_{H\otimes H} \colon k \to H\otimes H$ is the unit in monoid $(H \otimes H, \mu_{H \otimes H}, \eta_{H\otimes H})$ induced by $\eta_H$. 
		\item[(iii)] Equations  $$\epsilon\circ\eta_H = \eta_L,$$
		$$\epsilon\circ \mu_H\circ(\id_H\otimes(\alpha\circ\epsilon))
		= \epsilon\circ \mu_H = \epsilon\circ \mu_H\circ(\id_H\otimes(\beta\circ\epsilon))
		$$ 
	\end{itemize}
	hold. Morphism $\alpha$ is called \emph{source morphism} and $\beta$ is called \emph{target morphism}.
\end{definition}
This definition is good because of previously proven preparatory propositions.

In the category of vector spaces, condition (i) means that the image of the coproduct is inside Takeuchi product \cite{takeuchi}, and condition (ii) means that $\Delta(1_H) = 1_H \otimes_L 1_H$ and $\Delta(hh') = h_{(1)}h'_{(1)} \otimes_L h_{(2)}h'_{(2)}$, which is well defined because of (i). Condition (iii) in the category of vector spaces means that $\epsilon(1_H) = 1_L$ and $\epsilon(h\alpha(\epsilon(h'))) = \epsilon(h h') = \epsilon(h\beta(\epsilon(h')))$, or, equivalently, that the map defined by 
$h \btr l = \epsilon (h\alpha(l))$ and the map defined by $h \btr' l = \epsilon(h\beta(l))$ are left actions.

\begin{definition}
	Let $(R,\mu_R,\eta_R)$ be a monoid in $\mathcal{V}$. \emph{Internal right $R$-bialgebroid} in category $\mathcal{V}$ is given by the following data:
	\begin{itemize}
		\item[(i)] monoid $(H , \mu_H, \eta_H)$ in $\mathcal{V}$
		\item[(ii)] monoid morphisms  $\alpha\colon R \to H$, $\beta \colon R^{\op} \to H$ such that $$\mu_H \circ \tau_{H,H} \circ (\alpha \otimes \beta) = \mu_H \circ (\alpha \otimes \beta)$$
		\item[(iii)] comonoid $(H, \Delta, \epsilon)$ in the monoidal category of $R$-bimodules $({}_R \mathcal{V}{}_R, \otimes_R, R)$, where $H$ is an internal $R$-bimodule with regard to right action $\mu_H \circ (\id_H \otimes \alpha)$ and left coaction $\mu_H \circ \tau_{H,H} \circ (\beta \otimes \id_H)$. 
	\end{itemize}
	satisfying the following conditions:
	\begin{itemize}
		\item[(i)] Unique left action $\lambda$ for which the following diagram is commutative:
		$$
		\begin{tikzcd}
			(H \otimes H) \otimes (H \otimes H)  \ar{r}{\mu_{H\otimes H}} 
			\ar{d}{\id\otimes  \pi } &  
			(H\otimes H) 
			\ar{d}{ \pi} \\
			(H \otimes H) \otimes (H \otimes_R H)  \ar[dashed]{r}{\lambda} 
			& 
			(H\otimes_R H) 
		\end{tikzcd}
		$$
		satisfies $$\lambda \circ (\eta_H\otimes \beta \otimes\Delta)\circ(l_R^{-1}\otimes \id_H) = \lambda\circ (\alpha\otimes\eta_H \otimes \Delta)\circ(r_R^{-1}\otimes\id_H)
		,$$ where $r_R$ and $l_R$ are components of right unit and left unit of the monoidal category~$\mathcal{V}$, $\mu_{H\otimes H}$ is the multiplication in monoid $(H \otimes H, \mu_{H \otimes H}, \eta_{H\otimes H})$ induced by $\mu_H$ componentwise, and $\pi$ is the morphism in coequalizer $H\otimes H \to H\!\otimes_R\! H$.
		\item[(ii)] Unique morphism $\rho$ for which the following diagram is commutative: 
		$$
		\begin{tikzcd}
			(H\!\otimes\! H) \otimes H \ar{r}{\pi\otimes\id} \ar{d}{\id\otimes\Delta} & (H\!\otimes_R\! H) \otimes H \ar[dashed]{d}{\rho} \\
			(H\!\otimes\! H) \otimes (H\!\otimes_R\! H) \ar{r}{\lambda} & (H\!\otimes_R\! H)
		\end{tikzcd}
		$$
		is a right action. Equivalently,  $$\Delta\circ\eta_H = \pi\circ \eta_{H\otimes H}, $$ $$\Delta\circ \mu_H = \rho\circ (\Delta\otimes\id_H)$$ hold, where $\eta_{H\otimes H} \colon k \to H\otimes H$ is the unit in monoid $(H \otimes H, \mu_{H\otimes H}, \eta_{H\otimes H} )$ induced by $\eta_H$. 
		\item[(iii)] Equations $$\epsilon\circ\eta_H = \eta_R,$$
		$$\epsilon\circ \mu_H\circ( (\alpha\circ\epsilon) \otimes \id_H)= \epsilon\circ \mu_H = \epsilon\circ \mu_H\circ( (\beta\circ\epsilon) \otimes \id_H)$$
	\end{itemize}
	hold. Morphism $\alpha$ is called \emph{source morphism} and $\beta$ is called \emph{target morphism}. 
\end{definition}
This definition is good due to propositions for right bialgebroid that are analogous to preparatory propositions for the left bialgebroid. 

In the category of vector spaces, condition (i) means that the image of coproduct is inside Takeuchi product, and condition (ii) means that $\Delta(1_H) = 1_H \otimes_R 1_H$ and $\Delta(hh') = h_{(1)}h'_{(1)} \otimes_R h_{(2)}h'_{(2)}$, which is well defined because of (i). Condition (iii) in this category means that $\epsilon(1_H) = 1_R$ and $\epsilon(\alpha(\epsilon(h))h') = \epsilon(h h') = \epsilon(\beta(\epsilon(h))h')$, or, equivalently, that the map defined by
$r \btl h = \epsilon (\alpha(r)h)$ and the map defined by $r \btl' h = \epsilon(\beta(r)h)$ are right actions.

\section{Internal Hopf algebroid} \label{sec:internalHA}
\subsection{Preparatory results}
Let $R$ and $L$ be monoids in symmetric monoidal category $\mathcal{V}$ with coequalizers that commute with the monoidal product.

\begin{theorem}\label{lem:DeltaLRbim} Let $\mathcal{H}_L = (H, \mu_H,\eta_H, \alpha_L,\beta_L, \Delta_L, \epsilon_L)$ be an internal left  $L$-bialgebroid and let the structure of $R$-bimodule be given on $H$ by multiplication on the right with commuting monoid morphisms $\beta_R\colon R^{\op}\to H$, $\alpha_R \colon R\to H$. 
	If
	$$
	\alpha_L\circ\epsilon_L\circ\beta_R = \beta_R,
	\,\,\,\,\,\,\,\,\,
	\beta_L\circ\epsilon_L\circ\alpha_R = \alpha_R,
	$$
	then $\Delta_L$ is an $R$-bimodule morphism. 
	
	Analogously, for an internal right $R$-bialgebroid $\mathcal{H}_R = (H, \mu_H,\eta_H, \alpha_R,\beta_R, \Delta_R, \epsilon_R)$ with the structure of an $L$-bimodule on $H$ given by multiplying on the left with commuting monoid morphisms $\alpha_L \colon L\to H$, $\beta_L\colon L^{\op}\to H$ if
	$$
	\alpha_R\circ\epsilon_R\circ\beta_L = \beta_L,
	\,\,\,\,\,\,\,\,\,
	\beta_R\circ\epsilon_R\circ\alpha_L = \alpha_L,
	$$
	then $\Delta_R$ is an $L$-bimodule morphism. 
\end{theorem}

\begin{proof} We prove that $\Delta_L$ is a morphism of right $R$-modules; the proof that it is a morphism of left $R$-modules is analogous. Both claims are similarly proven for $\Delta_R$. 
	We prove that 
	$$\Delta_L\circ \mu_H \circ (\id_H\otimes\alpha_R) = 
	\nu_R \circ (\Delta_L \otimes \id_R),$$ where $\nu_R$ is the right $R$-action
	on $H\otimes_L H$ that is, by the definition of the monoidal product of a right $L$-module 
	and an $L$-$R$-bimodule, induced by right $R$-action $\mu_H \circ (\id_H\otimes\alpha_R)$ on $H$, that is, the unique morphism $\nu_R$ such that  
	$$\nu_R \circ (\pi_L \otimes \id_R) = \pi_L \circ (\id_H \otimes (\mu_H \circ (\id_H\otimes\alpha_R))).$$ 
	In the proof we imitate the following calculation done elementwise. First, we have that
	$$\Delta_L(h \cdot \alpha_R(r)) \stackrel{\text{(i)}}{=} "\Delta_L(h) \cdot \Delta_L(\alpha_R(r)) "$$
	and then we calculate
	$$\begin{array}{rl} 
	\Delta_L(\alpha_R(r)) & \stackrel{\text{(ii)}}{=} \Delta_L(\beta_L(\epsilon_L(\alpha_R(r))))
	\\ & \stackrel{\text{(iii)}}{=} \Delta_L(\beta_L(\epsilon_L(\alpha_R(r))) \cdot 1_H) 
	\\ & \stackrel{\text{(iv)}}{=}
	\Delta_L(1_H) \cdot \epsilon_L(\alpha_R(r)) 
	\\ & \stackrel{\text{(v)}}{=} (1_H \otimes_L 1_H) \cdot \epsilon_L(\alpha_R(r)) 
	\\ & \stackrel{\text{(vi)}}{=} 1_H \otimes_L (\beta_L(\epsilon_L(\alpha_R(r))) \cdot 1_H ) 
	\\ & \stackrel{\text{(vii)}}{=} 1_H \otimes_L \alpha_R(r)
	\end{array}$$
	and substitute this in the first equation to get
	$$"\Delta_L(h) \cdot \Delta_L(\alpha_R(r)) " = "\Delta_L(h) \cdot (1_H \otimes_L \alpha_R(r))" \stackrel{\text{(ix)}}{=} \Delta_L(h) \cdot \alpha_R(r).$$
	The two expressions inside quotes are actually  $$\lambda(h \otimes \Delta_L(\alpha_R(r))), $$ $$\lambda(h \otimes 1_H \otimes_L \alpha_R(r)) \stackrel{\text{(viii)}}{=} \rho(\Delta_L(h) \otimes 1_H \otimes \alpha_R(r)).$$
	Therefore, we use these facts in this order in the calculation that follows. 
	\begin{itemize}
		\item[(i)] $\lambda$ is a left action: $\Delta_L \circ \mu_H = \lambda \circ (\id_H \otimes \Delta_L ),$ 
		\item[(ii)] assumption  $\alpha_R = \beta_L \epsilon_L \alpha_R$ of the theorem, 
		\item[(iii)] unit axiom $\mu_H \circ (\id_H \otimes \eta_H) \circ r_H^{-1} = \id_H,$
		\item[(iv)] $\Delta_L$ is a morphism of right $L$-modules: $\Delta_L \circ  \mu_H \circ (\beta_L \otimes \id_H) \circ \tau_{H,L} = \nu_L \circ (\Delta_L \otimes \id_L),$ where $\nu_L \colon H \otimes_L H \otimes L \to H \otimes_L H$ is the right coaction induced by right coaction $\mu_H \circ (\beta_L \otimes \id_H) \circ \tau_{H,L} \colon H \otimes L \to H,$
		\item[(v)] $\lambda$ is a left action, hence for the unit we have $\Delta_L \circ \eta_H = \pi_L \circ \eta_{H \otimes H} = \pi_L \circ (\eta_{H} \otimes \eta_{H}) \circ l_k^{-1},$ 
		\item[(vi)] by the definition of right action $\nu_L$ and aforementioned right action, we have: 
		$\nu_L \circ (\pi_L \otimes \id_L) = \pi_L \circ (\id_H \otimes (\mu_H \circ (\beta_L \otimes \id_H) \circ \tau_{H,L})),$ 
		\item[(vii)] simple calculation below the main calculation, which uses the unit axiom $\mu_H \circ (\id_H \otimes \eta_H) = r_H,$
		\item[(viii)] by definition of $\lambda$ we have $\lambda \circ (\id_H \otimes \pi_L)=\rho \circ (\Delta_L \otimes \id_H \otimes \id_H),$
		\item[(ix)] additional calculation below the main calculation, which uses that $\nu_R$ is a unique morphism such that  $\nu_R \circ (\pi_L \otimes \id_R) = \pi_L \circ (\id_H \otimes (\mu_H \circ (\id_H\otimes\alpha_R)))$ and the fact that, by the definition of $\rho$, we have that $\rho \circ (\pi_L \otimes \id_{H\otimes H}) = \pi_L \circ \mu_{H\otimes H}$.
	\end{itemize}
	We have that
	$$\begin{array}{rcl}
	& & \Delta_L \circ \mu_H \circ (\id_H\otimes\alpha_R) 
	\\& \stackrel{\text{(i)}}{=} & \lambda \circ (\id_H \otimes (\Delta_L \circ \alpha_R))
	\\ &\stackrel{\text{(ii)}}{=}& \lambda \circ (\id_H \otimes (\Delta_L \circ \beta_L\epsilon_L\alpha_R))
	\\ &\stackrel{\text{}}{=}& \lambda \circ (\id_H \otimes (\Delta_L \circ \id_H \circ \beta_L\epsilon_L\alpha_R))
	\\ &\stackrel{\text{(iii)}}{=}& \lambda \circ (\id_H \otimes (\Delta_L \circ  \mu_H \circ (\id_H \otimes \eta_H) \circ r_H^{-1} \circ \beta_L\epsilon_L\alpha_R))
	\\ &\stackrel{\text{}}{=}& \lambda \circ (\id_H \otimes (\Delta_L \circ  \mu_H \circ (\id_H \otimes \eta_H) \circ (\beta_L\epsilon_L\alpha_R \otimes \id_k) \circ r_R^{-1} ))
	\\ &\stackrel{\text{}}{=}& \lambda \circ (\id_H \otimes (\Delta_L \circ  \mu_H \circ (\beta_L\epsilon_L\alpha_R \otimes \eta_H)\circ r_R^{-1} ))
	\\ &\stackrel{\text{}}{=}& \lambda \circ (\id_H \otimes (\Delta_L \circ  \mu_H \circ (\beta_L \otimes \id_H) \circ (\epsilon_L\alpha_R \otimes \eta_H)\circ r_R^{-1} ))
	\\ &\stackrel{\text{}}{=}& \lambda \circ (\id_H \otimes (\Delta_L \circ  \mu_H \circ (\beta_L \otimes \id_H) \circ \tau_{H,L} \circ ( \eta_H \otimes \epsilon_L\alpha_R) \circ \tau_{R,k} \circ r_R^{-1} ))
	\\ &\stackrel{\text{(iv)}}{=}& \lambda \circ (\id_H \otimes (\nu_L \circ (\Delta_L \otimes \id_L) \circ ( \eta_H \otimes \epsilon_L\alpha_R) \circ l_R^{-1} ))
	\\ &\stackrel{\text{}}{=}& \lambda \circ (\id_H \otimes (\nu_L \circ ((\Delta_L \circ \eta_H) \otimes \epsilon_L\alpha_R) \circ l_R^{-1} ))
	\\ &\stackrel{\text{(v)}}{=}& \lambda \circ (\id_H \otimes (\nu_L \circ ((\pi_L \circ \eta_{H\otimes H}) \otimes \epsilon_L\alpha_R) \circ l_R^{-1} ))
	\\ &\stackrel{\text{}}{=}& \lambda \circ (\id_H \otimes (\nu_L \circ (\pi_L \otimes \id_L) \circ (\eta_{H \otimes H} \otimes \epsilon_L\alpha_R)  \circ l_R^{-1} ))
	\\ &\stackrel{\text{}}{=}& \lambda \circ (\id_H \otimes (\nu_L \circ (\pi_L \otimes \id_L) \circ (\eta_{H} \otimes \eta_{H} \otimes \epsilon_L\alpha_R) \circ (l_k^{-1} \otimes \id_R) \circ l_R^{-1} ))
	\\ &\stackrel{\text{(vi)}}{=}& \lambda \circ (\id_H \otimes (\pi_L \circ (\id_H \otimes (\mu_H \circ (\beta_L \otimes \id_H) \circ \tau_{H,L})) \circ (\eta_{H} \otimes \eta_{H} \otimes \epsilon_L\alpha_R) \circ (l_k^{-1} \otimes \id_R) \circ l_R^{-1} ))
	\\ &\stackrel{\text{(vii)}}{=}& \lambda \circ (\id_H \otimes (\pi_L \circ (\eta_H \otimes (\alpha_R \circ l_R)) \circ (l_k^{-1} \otimes \id_R) \circ l_R^{-1} ))
	\\ &\stackrel{\text{}}{=}& \lambda \circ (\id_H \otimes (\pi_L \circ (\eta_H \otimes (\alpha_R \circ l_R)) \circ (\id_k \otimes l_R^{-1} ) \circ l_R^{-1} ))
	\\ &\stackrel{\text{}}{=}& \lambda \circ (\id_H \otimes (\pi_L \circ (\eta_H \otimes \alpha_R) \circ l_R^{-1} ))
	\\ &\stackrel{\text{}}{=}& \lambda \circ (\id_H \otimes \pi_L) \circ (\id_H \otimes ((\eta_H \otimes \alpha_R) \circ l_R^{-1} ))
	\\ &\stackrel{\text{(viii)}}{=}& \rho \circ (\Delta_L \otimes \id_H \otimes \id_H) \circ (\id_H \otimes ((\eta_H\otimes\alpha_R) \circ l_R^{-1}))
	\\ &\stackrel{}{=}& \rho \circ (\id_{H\otimes_L H} \otimes ((\eta_H\otimes\alpha_R) \circ l_R^{-1})) \circ (\Delta_H \otimes \id_R)
	\\ &\stackrel{\text{(ix)}}{=}& \nu_R \circ (\Delta_L \otimes \id_R).
	\end{array}
	$$
	The simple calculation for step (vii) is here:
	$$
	\begin{array}{rcl}
	&& (\id_H \otimes (\mu_H \circ (\beta_L \otimes \id_H) \circ \tau_{H,L})) \circ (\eta_{H} \otimes \eta_{H} \otimes \epsilon_L\alpha_R)  
	\\&=&\eta_H \otimes (\mu_H \circ (\beta_L \otimes \id_H) \circ \tau_{H,L} \circ (\eta_H \otimes \epsilon_L\alpha_R))
	\\&=&\eta_H \otimes (\mu_H \circ (\beta_L \otimes \id_H) \circ (\epsilon_L\alpha_R \otimes \eta_H) \circ \tau_{k,R})
	\\&=&\eta_H \otimes (\mu_H \circ (\beta_L \epsilon_L\alpha_R \otimes \eta_H) \circ \tau_{k,R})
	\\&=&\eta_H \otimes (\mu_H \circ (\alpha_R \otimes \eta_H) \circ \tau_{k,R})
	\\&=&\eta_H \otimes (\mu_H \circ (\id_H \otimes \eta_H) \circ (\alpha_R \otimes \id_k) \circ \tau_{k,R})
	\\&=&\eta_H \otimes (r_H \circ (\alpha_R \otimes \id_k) \circ \tau_{k,R})
	\\&=&\eta_H \otimes (\alpha_R \circ r_R \circ \tau_{k,R})
	\\&=&\eta_H \otimes (\alpha_R \circ l_R).
	\end{array}
	$$
	Finally, justification for step (ix) follows here. We prove that  
	$$\rho \circ (\id_{H\otimes_L H} \otimes ((\eta_H \otimes \alpha_R ) \circ l_R^{-1})) = \nu_R.$$
	Morphism $\nu_R$ is a unique morphism such that 
	$$\nu_R \circ (\pi_L \otimes \id_R) = \pi_L \circ (\id_H \otimes (\mu_H \circ (\id_H\otimes\alpha_R))).$$
	Let us prove that $\rho \circ (\id_{H\otimes_L H} \otimes ((\eta_H \otimes \alpha_R ) \circ l_R^{-1}))$ satisfies the same. We have that
	$$
	\begin{array}{rcl}
	&& \rho \circ (\id_{H\otimes_L H} \otimes ((\eta_H \otimes \alpha_R ) \circ l_R^{-1})) \circ   (\pi_L \otimes \id_R)
	\\&=&\rho \circ (\pi_L \otimes \id_{H\otimes H}) \circ (\id_{H\otimes H} \otimes \eta_H \otimes \alpha_R ) \circ (\id_{H\otimes H} \otimes l_R^{-1}) 
	\\&=&\pi_L \circ \mu_{H\otimes H} \circ (\id_{H\otimes H} \otimes \eta_H \otimes \alpha_R ) \circ (\id_{H\otimes H} \otimes l_R^{-1})
	\\&=&\pi_L \circ (\mu_{H} \otimes \mu_H) \circ (\id_{H} \otimes  \eta_H \otimes \id_{H} \otimes \alpha_R) \circ (\id_H \otimes \tau_{H,k} \otimes \id_R) \circ (\id_{H} \otimes \id_H \otimes l_R^{-1})
	\\&=&\pi_L \circ (r_H \otimes (\mu_{H} \circ ( \id_{H} \otimes \alpha_R)) \circ (\id_H \otimes \tau_{H,k} \otimes \id_R) \circ (\id_{H} \otimes \id_H \otimes l_R^{-1})
	\\&=&\pi_L \circ (\id_H \otimes (\mu_{H} \circ ( \id_{H} \otimes \alpha_R))) \circ (r_H \otimes \id_H \otimes \id_R) \circ (\id_H \otimes \tau_{H,k} \otimes \id_R) \circ (\id_{H} \otimes \id_H \otimes l_R^{-1}).
	\end{array}
	$$
	It remains to prove this claim involving unitors: 
	$$(r_H \otimes \id_H \otimes \id_R) \circ (\id_H \otimes \tau_{H,k} \otimes \id_R) \circ (\id_{H} \otimes \id_H \otimes l_R^{-1}) = \id_{H} \otimes \id_H \otimes \id_{R},$$
	but this easily follows from the axioms of symmetric monoidal category.
\end{proof}

\subsection{Definition of internal Hopf algebroid}
\begin{definition} 
	\emph{Internal Hopf algebroid} over base monoids $L, R$ in category $\mathcal{V}$ is given by the following data: internal left $L$-bialgebroid $\mathcal{H}_L = (H, \mu_H,\eta_H, \alpha_L,\beta_L, \Delta_L, \epsilon_L)$, internal right $R$-bialgebroid $\mathcal{H}_R =(H, \mu_H,\eta_H, \alpha_R,\beta_R ,\Delta_R, \epsilon_R)$ and monoid antihomomorphism $\tau \colon H \to H$ (antipode) which satisfy the following conditions. 
	\begin{itemize}
		\item[(i)]
		\begin{equation}\label{eq:alphaepsilonbeta}
		\begin{array}{lr}
		\alpha_L\circ\epsilon_L\circ\beta_R = \beta_R, \
		&
		\beta_L\circ\epsilon_L\circ\alpha_R = \alpha_R, 
		\\
		\alpha_R\circ\epsilon_R\circ\beta_L = \beta_L, \
		&
		\beta_R\circ\epsilon_R\circ\alpha_L = \alpha_L,
		\end{array}\end{equation}
		\item[(ii)]
		\begin{equation}\label{eq:deltalr1}
		(\Delta_R\otimes_L\id_H)\circ\Delta_L = 
		(\id_H\otimes_R\Delta_L)\circ\Delta_R,
		\end{equation}
		\begin{equation}\label{eq:deltalr2}
		(\Delta_L\otimes_R\id_H)\circ\Delta_R = 
		(\id_H\otimes_R \Delta_R)\circ\Delta_L.
		\end{equation}
		Monoidal product $\Delta_L\otimes_R\id_H$ is well defined because $\Delta_L$ is an $R$-bimodule morphism by Theorem~\ref{lem:DeltaLRbim} and similar claim holds for other monoidal products appearing above.
		\item[(iii)]
		\begin{equation}\label{eq:antipodeMain}
		\begin{array}{c}
		\tau\circ\beta_L =\alpha_L,\,\,\,\,\,\,\,\,\,\,\, 
		\tau\circ\beta_R = \alpha_R
		\\
		\mu_{\otimes'_L}\circ(\tau\otimes\id_H)\circ\Delta_L = \alpha_R\circ\epsilon_R ,
		\\
		\mu_{\otimes'_R}\circ(\id_H\otimes \tau)\circ\Delta_R = \alpha_L\circ\epsilon_L,
		\end{array}
		\end{equation}
		where $\mu_{\otimes'_R} \colon H \otimes'_R H \to H$ is the multiplication in $R$-ring $(H,\mu_{\otimes'_R},\alpha_R)$ induced by $\mu_H$ and $\mu_{\otimes'_L} \colon H \otimes'_L H \to H$ is the multiplication in $L$-ring $(H,\mu_{\otimes'_L},\alpha_L)$ induced by $\mu_H$, and the compositions of morphisms are well defined because the following morphisms  $\tau \otimes \id_H \colon H \otimes_L H \to H \otimes'_L H$ and $\id_H \otimes \tau \colon H \otimes _R H \to H \otimes'_R H$ are well defined.
	\end{itemize}
\end{definition}
\begin{proposition}
	In every $L$-bialgebroid $\epsilon_L\circ\alpha_L = \id_L = \epsilon_L \circ\beta_L$ holds. Furthermore,  $\alpha_L\circ\epsilon_L$ and $\beta_L\circ\epsilon_L$ are idempotent.
\end{proposition}
\begin{proof}
	We prove that $\id_L = \epsilon_L \circ\alpha_L$ by imitating the following calculation that is done elementwise:
	$$\epsilon_L(\alpha_L(l)) = \epsilon_L(\alpha_L(l) \cdot 1_H) = l \cdot \epsilon_L(1_H)= l \cdot 1_L = l.$$ Therefore, we use (from left to right) (i) unit axiom $\mu_L \circ(\id_L\otimes\eta_L)\circ r_L^{-1} = \id_L$, (ii) compatibility of counit and unit
	$\epsilon_L\circ\eta_H = \eta_L$, (iii)
	condition that $\epsilon_L$ is a morphism of left $L$-modules, $\epsilon_L\circ\mu_H\circ(\alpha_L\otimes\id_H)= \mu_L\circ(\id_L \otimes \epsilon_L) \colon L\otimes H\to L$, and (iv) unit axiom $\mu_H\circ(\id_H\otimes\eta_H)\circ r_H^{-1} = \id_H$. We have that
	$$\begin{array}{lcl}
	\id_L &\stackrel{\text{(i)}}{=}& \mu_L\circ (\id_L\otimes\eta_L)\circ r_L^{-1}
	\\ &\stackrel{\text{(ii)}}{=}&  \mu_L\circ (\id_L \otimes (\epsilon_L \circ \eta_H))\circ r_L^{-1}
	\\ &=& \mu_L\circ (\id_L \otimes \epsilon_L) \circ (\id_L \otimes  \eta_H)\circ r_L^{-1}
	\\ &\stackrel{\text{(iii)}}{=}& \epsilon_L\circ\mu_H\circ(\alpha_L\otimes\id_H)  \circ ( \id_L \otimes  \eta_H)\circ r_L^{-1}
	\\ &=& \epsilon_L\circ\mu_H\circ(\alpha_L\otimes\eta_H) \circ r_L^{-1}
	\\ &=& \epsilon_L\circ\mu_H\circ(\id_H \otimes \eta_H) \circ (\alpha_L\otimes\id_k) \circ r_L^{-1}
	\\ &\stackrel{\text{(iv)}}{=}& \epsilon_L\circ r_H \circ (\alpha_L\otimes\id_k) \circ r_L^{-1}
	\\ &=& \epsilon_L\circ\alpha_L.
	\end{array}$$
	In the last step we use the naturality of right unitor $r$.
	
	To prove that $\id_L = \epsilon_L \circ\beta_L$, we imitate the following  calculation that is done elementwise: 
	$$\epsilon_L(\beta_L(l)) = \epsilon_L(\beta_L(l) \cdot 1_H) = \epsilon_L(1_H) \cdot l = 1_L \cdot l = l.$$ Therefore, we use (from left to right) (i) unit axiom $\mu_L\circ(\eta_L\otimes\id_L)\circ l_L^{-1} = \id_L$, (ii) compatibility of counit and unit
	$\epsilon_L\circ\eta_H = \eta_L$, (iii)
	condition that $\epsilon_L$ is a morphism of right $L$-modules, $\epsilon_L\circ\mu_H\circ(\beta_L\otimes\id_H)\circ\tau_{H,L}= \mu_L\circ(\epsilon_L\otimes\id_L) \colon H\otimes L\to L$, and (iv) unit axiom $\mu_H\circ(\id_H\otimes\eta_H)\circ r_H^{-1} = \id_H$. We have that
	$$\begin{array}{lcl}
	\id_L &\stackrel{\text{(i)}}{=}& \mu_L\circ (\eta_L\otimes\id_L)\circ l_L^{-1}
	\\ &\stackrel{\text{(ii)}}{=}&  \mu_L\circ ((\epsilon_L \circ \eta_H)\otimes\id_L)\circ l_L^{-1}
	\\ &=& \mu_L\circ (\epsilon_L \otimes \id_L) \circ ( \eta_H \otimes\id_L)\circ l_L^{-1}
	\\ &\stackrel{\text{(iii)}}{=}& \epsilon_L\circ\mu_H\circ(\beta_L\otimes\id_H)\circ\tau_{H,L}  \circ ( \eta_H \otimes\id_L)\circ l_L^{-1}
	\\ &=& \epsilon_L\circ\mu_H\circ(\beta_L\otimes\id_H)\circ ( \id_L \otimes \eta_H) \circ \tau_{k,L}  \circ l_L^{-1}
	\\ &=& \epsilon_L\circ\mu_H\circ(\beta_L\otimes\id_H)\circ ( \id_L \otimes \eta_H) \circ r_L^{-1}
	\\ &=& \epsilon_L\circ\mu_H\circ(\beta_L\otimes\eta_H) \circ r_L^{-1}
	\\ &=& \epsilon_L\circ\mu_H\circ(\id_H \otimes \eta_H) \circ (\beta_L\otimes\id_k) \circ r_L^{-1}
	\\ &\stackrel{\text{(iv)}}{=}& \epsilon_L\circ r_H \circ (\beta_L\otimes\id_k) \circ r_L^{-1}
	\\ &=& \epsilon_L\circ\beta_L.
	\end{array}$$

\end{proof}

\begin{corollary}
	$\epsilon_R\circ\beta_L \colon L^{\op} \cong R$ is an isomorphism with inverse 
	$\epsilon_L\circ\alpha_R$.
\end{corollary}
\begin{proof}
	$\epsilon_L\alpha_R\epsilon_R\beta_L = \epsilon_L\beta_L=\id_L$ and
	$\epsilon_R\beta_L\epsilon_L\alpha_R = \epsilon_R\alpha_R = \id_R$.
\end{proof}

\begin{corollary}
	$\epsilon_L\circ\beta_R \colon R^{\op} \cong L$ is an isomorphism with inverse 
	$\epsilon_R\circ\alpha_L$.
\end{corollary}

\begin{proof}
	$\epsilon_R\alpha_L\epsilon_L\beta_R = \epsilon_R\beta_R=\id_R$ and 
	$\epsilon_L\beta_R\epsilon_R\alpha_L = \epsilon_L\alpha_L = \id_L$.
\end{proof}

\section*{Acknowledgements} The author is grateful to Zoran \v{S}koda for suggesting this approach and for discussions. Material in this article is translation of Chapter 7 of the author's dissertation of title \emph{Completed Hopf algebroids} \cite{CHA}.

\appendix

\section{Summary of \emph{Completed Hopf algebroids}} \label{summary} 
In the doctoral thesis \cite{CHA} of the author, a natural generalization of the definition of a Hopf algebroid is introduced, internal to any symmetric monoidal category with coequalizers that commute with the monoidal product. Furthermore, a symmetric monoidal category $(\ipV, \totimes, k)$ of filtered cofiltered vector spaces is constructed, 
whose morphisms are linear maps which in a weak sense respect 
the filtrations and cofiltrations, and whose monoidal product is 
the usual tensor product of vector spaces formally completed
and with a corresponding filtration of cofiltrations. It is proven then that this category satisfies the above conditions for the existence of internal Hopf algebroids. It contains two dual subcategories, the category $(\iV, \otimes, k)$ of filtered vector spaces and the category $(\pV,\hotimes, k)$ of cofiltered vector spaces. The monoidal product in it combines the ordinary tensor product and a completed tensor product.

An important class of Hopf algebroids over a noncommutative base
is comprised of smash products of a Hopf algebra $H$ and a braided commutative
algebra $A$ in the category of Yetter--Drinfeld modules over $H$. Such Hopf algebroids are called scalar extensions. In the thesis, it is proven that the smash products in which $H$ and $A$ are replaced by their analogues in the monoidal category of filtered cofiltered vector spaces have the structure of Hopf algebroids in that monoidal category. This sets the base for studying the Heisenberg doubles $A^* \sharp A$ in which $A$ is an infinite-dimensional Hopf algebra instead of a finite-dimensional one, among other examples, and the existence of the Hopf algebroid structure on them internal to the category $\ipV$.

Then Hopf pairings of a filtered Hopf algebra $A$ and a cofiltered Hopf algebra $H$ which are non-degenerate in the variable in $H$ are studied, and sufficient conditions for $A$ to be a braided commutative Yetter--Drinfeld module algebra over $H$ in the $\ipV$ category are found.  A smaller class of examples is also studied, for which $A$ is a Hopf algebra countably filtered by finite-dimensional vector spaces. Necessary and sufficient conditions on Hopf algebras  $A$ and $A^*$, or $A$ and $H$, are here found, in the form of finite dimensionality of the adjoint orbits of $A$ and the existence of certain canonical elements in $H \sharp A$. Thus a construction of some filtered cofiltered Hopf algebroids of scalar extension type is obtained. 
Important examples of such scalar extensions are the ones with $A$ the universal enveloping algebra $U(\gg)$ of a finite-dimensional Lie algebra $\gg$.   
When $H$ is equal to its algebraic dual $\Ug^*$ with induced cofiltration, the corresponding scalar extension, that is the Heisenberg double of $\Ug$, can be identified as an algebra with the algebra
of differential operators on the formal neighborhood of the unit
of a Lie group integrating $\gg$, suggesting applications in geometry
and mathematical physics.

\bigskip\bigskip
\noindent Martina Stoji\'{c} \\
\emph{Address:} Department of Mathematics, Faculty of Science, University of Zagreb, Bijeni\v{c}ka cesta 30, 10000 Zagreb, Croatia \\
\emph{E-mail address:} \texttt{stojic@math.hr} \\
\emph{ORCID:} \url{https://orcid.org/0000-0002-7994-7509}
\end{document}